\documentclass[a4paper,12pt]{article} 
\usepackage{amsmath, amsthm, amssymb}
\usepackage{url}

\usepackage{physics}
\usepackage{graphicx,lipsum}
\graphicspath{ {./downloads/} }

\usepackage[colorlinks,citecolor=red,urlcolor=blue,bookmarks=false,hypertexnames=true]{hyperref}
\usepackage{tikz}
\usetikzlibrary{calc}
\usetikzlibrary{shapes}
\usepackage[autostyle]{csquotes}
\makeatletter

\usepackage[colorlinks]{hyperref}
\usepackage[nameinlink,capitalize]{cleveref}
\newtheorem{theorem}{Theorem}[section]
\newtheorem{corollary}[theorem]{Corollary}
\newtheorem{lemma}[theorem]{Lemma}
 \newtheorem{proposition}[theorem]{Proposition}
 \newtheorem{remark}[theorem]{Remark}
\newtheoremstyle{named}{}{}{\itshape}{}{\bfseries}{.}{.5em}{\thmnote{#3's }#1}
\theoremstyle{named}

\theoremstyle{definition}
\newtheorem{definition}[theorem]{Definition}
\newtheorem{example}[theorem]{Example}

\usepackage{xspace}
\usepackage[margin=1.0in]{geometry}

\usepackage{titlesec}

\usepackage{mathtools}

\DeclarePairedDelimiterX{\inp}[2]{\langle}{\rangle}{#1, #2}
\titleformat{\chapter}
  {\Large\bfseries} 
  {}                
  {0pt}            
  {\huge}

\begin{document}

\begin{center}
\fontsize{13pt}{10pt}\selectfont
    \textsc{\textbf{ ON GRADED CLASSICAL S-PRIMARY SUBMODULES}}
    \end{center}
\vspace{0.1cm}
\begin{center}
   \fontsize{12pt}{10pt}\selectfont
    \textsc{{\footnotesize Tamem \ Al-Shorman and\  Malik\ Bataineh \  }}
\end{center}
\vspace{0.2cm}

\begin{abstract}
  The purpose of this article is to introduce the graded classical S-primary submodules which are extensions of graded classical primary submodules. We state that P is a graded classical S-primary submodule of R-module M if there exists  $s\in S$ such that $x,y \in h(R)$ and $m \in h(M)$, if $xym \in P$, then $sxm \in P$ or $sy^nm \in P$ for some positive integer n. Several properties and characteristics of graded classical S-primary submodules have been studied.    
\end{abstract}

\section{Introduction}
Through this article, we assume that R is a commutative graded ring with nonzero unity and M is a unitary graded R-module. 
 
Let G be an abelian group with identity e. A ring R is called a G-graded ring if $ R= \bigoplus\limits_{g \in G} R_g$   with the property $R_gR_h\subseteq R_{gh}$ for all $g,h \in G$, where $R_g$ is an additive subgroup of R for all $g\in G$. The elements of $R_g$ are called homogeneous of degree g. If $x\in R$, then $x$ can be written uniquely as $\sum\limits_{g\in G} x_g$, where $x_g$ is the component of $x$ in $R_g$. The set of all homogeneous elements of R is $h(R)= \bigcup\limits_{g\in G} R_g$. Let P be an ideal of a G-graded ring R. Then P is called a graded ideal if $P=\bigoplus\limits_{g\in G}P_g$, i.e, for $x\in P$ and  $x=\sum\limits_{g\in G} x_g$ where $x_g \in P_g$ for all $g\in G$. An ideal of a G-graded ring is not necessary graded ideal (see \cite{abu2019graded}).  

Let R be a G-graded ring . A left R-module M is said to be a graded R-module if there exists  a family of additive subgroups $\{M_g\}_{g \in G}$ of M such that $M= \bigoplus\limits_{g \in G} M_g$ with the property $M_gM_h\subseteq M_{gh}$ for all $g,h \in G$. The set of all homogeneous elements of M is $h(M)= \bigcup\limits_{g\in G} M_g$. Note that $M_g$ is an $R_e$-module for every $g\in G$. A submodule P of M is said to be graded submodule of M if $P=\bigoplus\limits_{g\in G}P_g$. In this case, $P_g$ is called g-component of P (see \cite{gordon1982graded}).

Let P be a proper graded ideal of R. Then the graded radical of P is denoted by Grad(P) and it is defined as written below:
\begin{center}
  {\small $Grad(P)=\Big{\{} x= \sum\limits_{g\in G} x_g \in R$ : for all $g\in G$, there exists s $n_g\in \mathbb{N}$ such that ${x_g}^{n_g} \in P \Big{\}}$}.  
\end{center}
Note that Grad(P) is always a graded ideal of R (see
\cite{refai2000graded}). Let P be a graded submodule of M and I is a graded ideal of R. Then $(P:_RM)$ is defined as $(P:_RM) = \{r\in R: rM\subseteq P \}$ and then $(P:_RM)$ is a graded submodule of M (see \cite{atani2006graded}). The graded submodule $(P:_MI)$ is defined as $(P:_MI)=\{m\in M : Im \subseteq P\}$. Particularly, we use $(P:_Ms)$ instead of $(P:_MRs)$. 

Refai and Al-Zoubi in \cite{refai2004graded} introduced the notion of graded primary ideals. The definition of graded primary submodules has piqued the interest of a number of mathematicians, see for example \cite{al2016graded} and \cite{atani2007graded}. Then many generalizations of graded primary submodule were investigated, such as graded classical primary submodule (see \cite{al2017graded}). A graded submodule P of M is said to be graded classical primary submodule of M if $p \neq M$, and whenever $x,y \in h(R)$ and $m \in h(M)$ with $xym \in P$, then either $xm \in P$ or $y^nm \in P$ for some positive integer n.

Let $S\subseteq h(R)$ be a multiplicative closed set of R and P is a graded submodule of M such that $(P:_RM) \cap S = \emptyset$. Recently, Al-Zoubi and Ababneh in \cite{al2021graded}, defined the concept of graded S-primary submodules which is another generalization of graded primary submodules. A graded submodule P of M is said to be a graded S-primary submodule of M if there exists  $s\in S$ and whenever $xm\in P$, whree $x \in h(R)$ and $m \in h(M)$, then $sx \in Grad(P:_RM)$ or $sm \in P$. Also, we study the relationship between graded classical S-primary submodule and graded classical S-prime submodule. In Section Two, we introduce the concept of graded classical S-primary submodule. We say that P is a graded classical S-primary submodule of R-module M if there exists  $s\in S$ such that $x,y \in h(R)$ and $m \in h(M)$, if $xym \in P$, then $sxm \in P$ or $sy^nm \in P$ for some positive integer n. 

\section{Graded Classical S-Primary Submodules}
In this section, we introduce the notation of graded classical S-primary submodules and show that graded classical S-primary submodules have a lot of similar features to these of graded classical primary submodules. 
\begin{definition}\label{def1}
Let $S\subseteq h(R)$ be a multiplucative closed set of R and let P be a graded ideal of R such that $P \cap S =\emptyset$. A graded ideal P of R is said to be a graded classical S-primary ideal of R if there exists  $s\in S$ and whenever $xy \in P$, where $x,y \in h(R)$ then $sx \in P$ or $sy^n  \in P$ for some positive integer n.
\end{definition}

\begin{definition}\label{def2}
Let $S \subseteq h(R)$ be a multiplicative closed set of R and let P be a graded submodule of M such that $(P:_R M) \cap S = \emptyset$ . A graded submodule P of M is said to be a graded classical S-primary submodule of M if there exists  $s\in S$ and whenever $xym\in P$, where $x,y \in h(R)$ and $m \in h(M)$ then $sxm \in P$ or $sy^n m \in P$ for some positive integer n.   
\end{definition}

\begin{lemma}\label{lem1}
Let M be a graded R-module and $S\subseteq h(R)$ be a multiplicative  closed set of R. Then the following is true.
\\
1. If P is a graded classical primary submodule of M such that $(P:_R M) \cap S = \emptyset$, then P is a graded classical S-primary submodule of M.
\\
2. If P is a graded classical S-primary submodule of M and $S \subseteq U(R)$, where U(R) denoted the set of units in R, then P is a graded classical primary submodule of M.
\\
\\
\textbf{Proof.} This is clear by using the definitions of graded classical primary submodules and graded classical S-primary submodules.
\end{lemma}

Clearly, every graded S-primary submodule is a graded classical S-primary submodule. However, the converse is not true. 

\begin{example}\label{exm1}
Let R be a graded integral domain, $S\subseteq h(R)$ be a multiplicative closed set of R and I be a non-zero graded prime ideal of R with $I \cap S = \emptyset$. In this case $P = I \oplus (0) $ is a graded classical primary submodule of the graded R-module $R\oplus R$ so P is a graded classical S-primary submodule of the graded R-module  while it is not graded S-primary submodule.
\end{example}

Let $S\subseteq h(R)$ be a multiplicative closed set of R and P is a graded submodule of R-module M with $(P:_RM)\cap S = \emptyset$. Then P is said to be graded classical S-prime submodule of M if there exists  $s\in S$ and whenever $xym \in P$ where $x,y \in h(R) $ and $m\in h(M)$, then $sxm \in P$ or $sym \in P$ (see \cite{kmm2021graded}).

\begin{remark}\label{rem1}
Every graded classical S-prime submodule is a graded classical S-primary submodule. However, the next theorem shows that graded classical S-primary submodule need not to be graded classical S-prime submodule. This can be seen by taking $S = U(R) \cap h(R)$ and recalling that graded classical primary submodule does not imply graded classical prime submodule.
\end{remark}

\begin{theorem}\label{thm1}
Let $S\subseteq h(R)$ be a multiplicative closed set of R. If P is a graded classical S-primary submodule of R-module M, then $S^{-1}P$ is a graded classical primary submodule of $S^{-1}R$-module $S^{-1}M$. 
\\
\\
\textbf{Proof.} Let $(\frac{x}{s_1})(\frac{y}{s_2})(\frac{m}{s_3})\in S^{-1}P$, where $x, y \in h(R)$,  $ m\in h(M)$ and $s_1, s_2, s_3 \in S$. Then $xym \in P$ since P is a graded classical S-primary submodule of M, so there exists  $s\in S$ such that $sxm \in P $ or $sy^nm \in P$ for some positive integer n. This yields that $\frac{x}{s_1} \frac{m}{s_3} = \frac{s x m}{ss_1s_3} \in S^{-1}P$ or $\frac{y^n}{s_2}\frac{m}{s_3} = \frac{sy^nm}{ss_2s_3} \in S^{-1}P$ for some positive integer n. Hence, $S^{-1}P$ is a graded classical primary submodule of $S^{-1}M$. 
\end{theorem}

\begin{lemma}\label{lem2}
Let $S_1 \subseteq S_2$ be a multiplucative subset of h(R) and P is a graded submodule of R-module M such that $(P:_R M) \cap S_2 =\emptyset$. If P is a graded classical $S_1$-primary submodule of M, then P is a graded classical $S_2$-primary submodule of M.
\\
\\
\textbf{Proof.} Let $xym \in P$, where $x, y \in R$ and $m\in h(M)$. Since P is a graded classical $S_1$-primary submodule of M, so there exists  $s \in S_1$ such that $sxm \in P$ or $sy^nm \in P$ for some positive integer n. But $S_1 \subseteq S_2$ then $s\in S_2$. Hence, P is a graded classical $S_2$-primary submodule of M.  
\end{lemma}

\begin{proposition}\label{prop1}
Let $S_1, S_2 \subseteq h(R)$ be a multiplicative set of R with $S_1 \subseteq S_2$ and P be a graded submodule of R-module M such that $(P:_RM) \cap S_2 = \emptyset$. For any $s\in S_2$, there is $r\in S_2$ satisfying $sr\in S_1$. If P is a graded classical $S_2$-primary submodule of M, then P is a graded classical $S_1$-primary submodule of M.
\\
\\
\textbf{Proof.} Let $xym \in P $, where $x, y \in h(R)$ and $m\in h(M)$. Since P is a graded classical $S_2$-primary submodule of M, so there exists  $s\in S_2$ then $sxm \in P$ or $sy^nm \in P$ for some positive integer n. Assume that $k = sr \in S_1$ for some $r \in S_2$, then $kxm\in P$ or $ky^nm \in P$. Hence, P is a graded classical $S_1$-primary submodule of M 
\end{proposition}

\begin{remark}\label{rem2}
Let S be a multiplicative closed subset of a graded ring R and $S^*$ denoted the saturation of S is defined by 
\begin{center}
    $S^* =\{r \in h(R): \frac{r}{1}$ is unit in $S^{-1}R \}$
\end{center}

It is clearly that $S^*$ is a multiplicative closed subset of a graded ring R containing S ( see \cite{gilmer1992multiplicative}).
\end{remark}

\begin{proposition}\label{prop2}
Let $S\subseteq h(R)$ be a multiplicative closed set of R and P is a graded submodule of R-module m such that $(P:_RM) \cap S = \emptyset$. Then P is a graded classical S-primary submodule of R-module M if and only if P is a graded classical $S^*$-primary submodule of M.
\\
\\
\textbf{Proof.} ($\Rightarrow$) Clearly, since $S \subseteq S^*$ and $(P:_RM) \cap S^* =\emptyset$. By Lemma \ref{lem2}, P is a graded classical $S^*$-primary submodule of M.
\\
($\Leftarrow$) It follows by Proposition 2.18 in \cite{al2022generalizations} and Proposition \ref{prop1}.
\end{proposition}

\begin{theorem}\label{thm2}
Let $S\subseteq h(R)$ be a multiplicative closed set of R and $P_i$ for all $i =\{1, ..., n\}$ be a graded submodules of R-module M such that $(P:_RM) \cap S = \emptyset$. If $P_i$ is a graded classical S-primary submodule of M for each i, then $\bigcap\limits_{i=1}^{n} P_i$ is a graded classical S-primary submodule of M.
\\
\\
\textbf{Proof.} Let $xym \in P_i$ for each i, where $x, y \in h(R)$ and $m \in h(M)$. Since $P_i$ is a graded classical S-primary submodule  of M, so there exists  $s_i \in S$ such that $s_ixm \in P_i$ or $s_iy^nm \in P_i$ for some positive integer n. Now take $s = \prod\limits_{i=1}^{n} s_i$ then $s\in S$. Suppose that $x, y \in h(R)$, $m \in h(M)$ and $xym \in \bigcap\limits_{i=1}^{n} P_i$, then $xym \in P_i$ for all $i=\{ 1, ..., n\}$. But $P_i$ is a graded classical S-primary submodule of M, so $sxm \in P_i$ or $sy^nm \in P_i$ for each i. Thus $sxm \in \bigcap\limits_{i=1}^{n} P_i $ or $sy^nm \in \bigcap\limits_{i=1}^{n} P_i$. Hence, $\bigcap\limits_{i=1}^{n} P_i$ is a graded classical S-primary submodule of M. 
\end{theorem}

Recall that, if M and T be two G-graded R-modules. Then an R-module homomorphism $f : M \rightarrow T$ is said to be a G-graded R-module homomorphism if $f(M) \subseteq T$ (see \cite{nastasescu2004methods}).

\begin{theorem}\label{thm3}
Let $S \subseteq h(R)$ be a multiplicative closed set of R, $M_1, M_2$ be two graded R-module and $f: M_1 \rightarrow M_2$ be a graded R-module homomorphism. Then the following hold:
\\
(i) If f is a graded epimorphism and P is a graded classical S-primary submodule of $M_1$ with $ker(f) \subseteq P$, then f(P) is a graded classical S-primary submodule of $M_2$.
\\
(ii) If f is a graded epimorphism and I is a graded classical S-primary submodule of $M_2$, then $f^{-1}(I)$ is a graded classical S-primary submodule of $M_1$.
\\
\\
\textbf{Proof.} (i) At first, we need to show that $(f(P):_R M_2) \cap S = \emptyset$. Let $r\in (f(P):_R M_2) \cap S $, so $r M_2 \subseteq f(P)$ and then $f(r M_1)= rf(M_1)=r M_2 \subseteq f(P)$. Then $r M_1 = r M_1 +ker(f) \subseteq f^{-1}(f(P))= P + ker(f) = P$. Thus $r \in (P:_R M_1)$, which is contradiction since 
$(P:_RM_1) \cap S = \emptyset$. Now, let $xym_2 \in f(P)$, where $x, y \in h(R)$ and $m_2 \in h(M_2)$. Then there some $m_1 \in h(M_1)$ such that $m_2 = f(m_1)$. Thwn $xy m_2 = xy f(m_1)= f(xym_1) \in f(P)$. Thus there exists  $p \in P \cap h(M_1)$ such that $f(xym_1) = f(p)$. Hence $xym_1 - p \in ker(f) \subseteq P$ and then $xym_1 \in P$. Since P is a graded classical S-primary submodule of $M_1$, then there exists  $s \in S$ such that $sxm_1 \in P$ or $sy^nm_1 \in P$ for some positive integer n. So, $f(sxm_1) = sx f(m_1)=sxm_2 \in f(P)$ or $f(sy^nm_1)=sy^nf(m_1)=sy^nm_2 \in f(P)$. Therefore f(P) is a graded classical S-primary submodule of $M_2$.
\\
(ii) First, we need to show that $(f^{-1}(I):_R M_1) \cap S =\emptyset$. Let $r\in (f^{-1}(I):_R M_1) \cap S $, so $rM_1 \subseteq f^{-1}(I)$ and so $rf(M_1)= rM_2 \subseteq I$, this is a contradiction since $(I:_RM_2) \cap S =\emptyset$. Now, let $xy m_1 \in f^{-1}(I)$, where $x,y \in h(R)$ and $m_1 \in h(M_1)$. Thus $f(xym_1)= xy f(m_1)\in I$. But I is a graded classical S-primary submodule of $M_2$ then there exists  $s\in S$ such that $sxf(m_1) \in I$ of $sy^n f(m_1) \in I$ for some positive integer n. Hence $sxm_1 \in f^{-1}(I)$ or $sy^nm_1 \in f^{-1}(I)$. Therefore $f^{-1}(I)$ is a graded classical S-primary submodule of $M_1$.
\end{theorem}

\begin{proposition}\label{prop3}
Let P be a graded submodule of M and $S\subseteq h(R)$ be a multiplucative closed set of R. If P is a graded classical S-primary submodule of M, then $(P:_RM)$ is a graded classical S-primary ideal of R.
\\
\\
\textbf{Proof.} Let $xy\in (P:_RM)$, where $x,y \in h(R)$. Then $xym \in P$ for all $m \in h(M)$. Since P is a graded classical S-primary submodule of M, so there exists  $s\in S$ such that $sxm \in P $ or $sy^nm \in P$ for some positive integer n. Hence $sx \in (P:_RM)$ or $sy^n \in (P:_RM)$. Therefore $(P:_RM)$ is a graded classical S-primary ideal of R. 
\end{proposition}

\begin{proposition}\label{prop4}
Let $R = R_1\times R_2$ be a graded ring and $S = S_1 \times S_2 \subseteq h(R)$ be a multiplicative closed set of R. Assume that $P = P_1 \times P_2$ is a graded ideal of R. Them the following are equivalent:
\\
(i) P is a graded classical S-primary ideal of R.
\\
(ii) $P_1$ is a graded classical $S_1$-primary ideal of $R_1$ and $P_2 \cap S_2 \neq \emptyset$ or $P_2$ is a graded classical $S_2$-primary ideal of $R_2$ and $P_1 \cap S_1 \neq \emptyset$.
\\
\\
\textbf{Proof.} (i) $\Rightarrow$ (ii): Since $(a,0)(0,b) =(0,0) \in P$, so there exists  $s=(s_1,s_2) \in S$ and so $s(a,0)= (s_1a,0)\in P$ or $s(0,b)^n = (0, s_2 b^n) \in P$ for some positive integer n, then $P_1 \cap S_1 \neq \emptyset$ or $P_2\cap S\neq \emptyset$. Suppose that $P_1 \cap S_1 \neq \emptyset$. As $P \cap S = \emptyset$ so $P_2 \cap S_2 = \emptyset$. Let $xy \in P_2$ for some $x,y \in h(R_2)$. Since $(0,x)(0,y) \in P$ and P is a graded classical S-primary ideal of R. We have $s(0,x)=(0,s_2x)\in P$ or $s(0,y)^n = (0,s_2y^n) \in P$ for some positive integer n. This yields $s_2x \in P_2$ or $s_2y^n \in P_2$. Therefore $P_2$ is a graded classical $S_2$-primary ideal of $R_2$. On the other hand, it is simple to prove that $P_1$ is a graded classical $S_1$-primary ideal of $R_1$.
\\
(ii) $\Rightarrow$ (i): Suppose that $P_1$ is a graded classical $S_1$-primary ideal of $R_1$ and $P_2 \cap S_2 \neq \emptyset$. Then there exists  $s_2 \in P_2 \cap S_2$. Let $(x,y)(a,b) = (xa,yb)\in P$ for some $x, a \in h(R_1)$ and $y,b \in h(R_2)$. Thus $xa \in P_1$ and thus there exists  $s_1 \in S_1$ then $s_1x \in P_1$ or $sa^n \in P_1$ for some positive integer n. Take $s=(s_1,s_2)\in S$ then take notice of this $s(x,y)=(s_1x,S_2y) \in P$ or $s(a,b)^n= (s_1a^n, s_2b^n) \in P$. Therefore P is a graded classical S-primary ideal of R. In other case, one can similarly prove that P is a graded classical S-primary ideal of R.
\end{proposition}

\begin{theorem}\label{thm4}
Let $S_1 , S_2 $ be a multiplicative closed subsets of rings $R_1 , R_2$ respectively and $P_1, P_2$ be a graded submodules of an $R_1$-module $M_1$ and $R_2$-module $M_2$ respectively. Suppose that $M = M_1 \times M_2$ as $R=R_1 \times R_2$-module, $S = S_1\times S_2$ and $P = P_1 \times P_2$. Then the following are equivalent:
\\
(i) P is a graded classical S-primary submodule of M.
\\
(ii) $P_1$ is a graded classical $S_1$-primary submodule of $M_1$ and $(P_2:_{R_2} M_2) \cap S_2 \neq \emptyset$ or $P_2$ is a graded classical $S_2$-primary submodule of $M_2$ and $(P_1:_{R_1} M_1) \cap S_1 \neq \emptyset$.
\\
\\
\textbf{Proof.} (i) $\Rightarrow$ (ii): By Proposition \ref{prop3}, $(P:_R M) = (P_1:_{R_1} M_1) \times (P_2:-{R_2} M_2)$ is a graded classical S-primary ideal of R and by Proposition \ref{prop4}, either $(P_1:_{R_1} M_1) \cap S_1 \neq \emptyset$ or $(P_2:_{R_2} M_2) \cap S_2 \neq \emptyset$. Suppose that  $(P_2:_{R_2} M_2) \cap S_2 \neq \emptyset$. Now, we will show that $P_1$ is a graded classical $S_1$-primary submodule of $M_1$. Let $xym_1 \in P_1$ for some $x,y \in h(R_1)$ and $m_1 \in h(M_1)$. Then $(x,1)(y,0)(m_1,0)=(xym,0) \in P$. As P is a graded classical S-primary submodule of M, there is an $s=(s_1,s_2) \in S$, so $s(x,1)(m_1,0)=(s_1xm_1,0)\in P$ or $s(y,0)^n(m_1,0) = (s_1y^nm_1,0)\in P$. This is implies that $s_1xm_1 \in P_1$ or $s_1y^nm_1 \in P_1$. Therefore $P_1$ is a graded classical $S_1$-primary submodule of $M_1$. On the other hand, it can be similarly show that  $P_2$ is a graded classical $S_2$-primary submodule of $M_2$.
\\
(ii) $\Rightarrow$ (i): Suppose that $(P_2:_{R_2} M_2) \cap S_2 \neq \emptyset$ and $P_1$ is a graded classical $S_1$-primary submodule of $M_1$. Then there exists  $ s_2 \in (P_2:_{R_2} M_2) \cap S_2$. Let $(x,y)(a,b)(m_1,m_2) = (xam_1 , ybm_2) \in P$, where $x,a \in h(R_1)$, $y,b \in h(R_2)$ and $m_i \in h(M_i)$ where $i = 1,2$. Then $xam_1 \in P_1$. As $P_1$ is a graded classical $S_1$-primary submodule of $M_1$ so there exists  $s_1 \in S_1$ and so $s_1xm_1 \in P_1$ or $s_1a^n m_1 \in P_1$. Now, take $s=(s_1,s_2) \in S$ then take notice of this $s(x,y)(m_1,m_2)= (s_1xm_1,s_2ym_2)\in P$ or $s(a,b)^n(m_1,m_2)=(s_1a^nm_1, s_2b^nm_2) \in P$. Therefore, P is a graded classical S-primary submodule of M. In other case, one can similarly prove that P is a graded classical S-primary submodule of M.  
\end{theorem}

\begin{theorem}\label{thm5}
Let $R=R_1 \times R_2 \times ... \times R_n$ where each $R_i$ is a graded ring for all $i=\{1,2,...,n\}$, $S=S_1\times S_2 \times ... \times S_n \subseteq h(R)$ be a multiplucative closed set of R and $M = M_1 \times M_2 \times ... \times M_n$ be R-module. Assume $P = P_1 \times P_2 \times ... \times P_n$ is a graded submodule of M. Then the following are equivalent:
\\
(i) P is a graded classical S-primary submodule of M.
\\
(ii) $P_i$ is a graded classical $S_i$-primary submodule of $M_i$ for some $i \in \{1,2,...,n\}$ and $(P_j:_{R_j} M_j) \cap S_j \neq \emptyset$ for all $j \in \{1,2,...,n\}-\{i\}$.
\\
\\
\textbf{Proof.} For n = 1, the result is true. If n = 2, then (i) $\Leftrightarrow$ (ii) is true by Theorem \ref{thm4}. Suppose that (i) and (ii) are equivalent when $k < n$. Now, we shall prove (i) $\Leftrightarrow$ (ii) when k = n. Let $P = P_1 \times P_2 \times ... \times P_n$. Take $I = P_1 \times P_2 \times ... \times P_{n-1}$ and $T = S_1\times S_2 \times ... \times S_{n-1}$ . Then by Theorem \ref{thm4}, a condition that is both necessary and sufficient for $P = I \times P_n$ is a graded classical S-primary submodule of M is that I is a graded classical S-primary submodule of $M_I$ and  $(P_n:_{R_n} M_n) \cap S_n \neq \emptyset$, where $M_I = M_1 \times M_2 \times ... \times M_{n-1} $ and $R_I = R_1 \times R_2 \times ... \times R_{n-1}$. The rest of the argument is based on the induction hypothesis.
\end{theorem}

\begin{theorem}\label{thm6}
Let $S \subseteq h(R)$ be a multiplucative closed set of R and P is a graded submodule of R-module M with $(P:_RM) \cap S = \emptyset$. Then the following are equivalent:
\\
(i) P is a graded classical S-primary submodule of M.
\\
(ii) there exists  $s\in S$ such that whenever $xI \subseteq P$ and I is a graded submodule of M and $x\in h(R)$, then $sx \in Grad(P:_R M) $ or $sI \subseteq P$.
\\
(iii) there exists  $s\in S$ such that whenever $JI \subseteq P$, where I is a graded submodule of M and J is a graded ideal of R, then $sJ \subseteq Grad(P:_RM)$  or $sI \subseteq P$.
\\
\\
\textbf{Proof.} (i) $\Rightarrow$ (ii): It is clear.
\\
(ii) $\Rightarrow$ (iii): Let $JI \subseteq P$, where I is a graded submodule of M and J is a graded ideal of R. We need to show that there exists  $s\in S$ such that $sI \subseteq P$ or $sJ \subseteq Grad(P_:RM)$. Since $xI \subseteq P$ for every $x\in J$, there exists  $s\in S$ such that$sx \in Grad(P:_RM) $ or $sI \subseteq P$ for every $x\in J$.
\\
(iii) $\Rightarrow$ (i): Let $x,y \in h(R)$ and $m \in h(M)$ with $xym \in P$. Now, set $J = Ry$ and $I = Rxm$, then I is a graded submodule of M and J is a graded ideal of R. Then we might get to the conclusion that $JI = Rxym \subseteq P$. then there exists  $s \in S$ so that $sJ = Rys \subseteq Grad(P:_R M) $ or $sI = Rxms \subseteq P $ and hence either $sy^nm \in P$ or $sxm \in P$. Therefore, P is a graded classical S-primary submodule of M.
\end{theorem}

Recall that, a graded R-module M over graded ring R is said to be a graded
multiplication module (gr-multiplication module) if for every graded submodule P of M there exists s a graded ideal I of R such that P = IM. It is clear that M is gr-multiplication R- module if and only if $P = (P :_R M)M$ for every graded submodule P of M (see \cite{escoriza1998multiplication}).

\begin{proposition}\label{prop5}
Let P be a graded submodule of M and $S\subseteq h(R)$ be a multiplucative closed set of R. If M is a graded multiplication module and $(P:_RM)$ is a graded classical S-primary ideal of R, then  P is a graded classical S-primary submodule of M.
\\
\\
\textbf{Proof.} Suppose that  M is a graded multiplication module and $(P:_RM)$ is a graded classical S-primary ideal of R. Let J be a graded ideal of R and I be a graded submodule of M with $JI \subseteq P$. Then $J(I:_RM) \subseteq (JI:_RM) \subseteq (P:_RM)$. As $(P:_RM)$ is a graded classical S-primary ideal of R, there is $s\in S $ so that $sI  = s(I:_RM)M \subseteq (P:_RM)M = P$ or $sJ \subseteq Grad(P:_RM)$. Therefore, by Theorem \ref{thm5} (iii), P is a graded classical S-primary submodule of M. 
\end{proposition}

\begin{lemma}\label{lem3}
Let $S\subseteq h(R)$ be a multiplicative closed set of R and P is a graded classical S-primary submodule of R-module M. Then the following hold for some $s\in S$.
\\
(i) $(P:_M s') \subseteq (P:_Ms)$ for all $s' \in S$.
\\
(ii) $((P:_RM):_Rs') \subseteq ((P:_RM):_Rs)$ for all $s' \in S$.
\\
\\
\textbf{Proof.} Similar to Lemma 2.1 in \cite{al2021graded}.  
\end{lemma}

\begin{theorem}\label{thm7}
Let M be a finitely generated graded R-module, $S \subseteq h(R)$ be a multiplicative closed set of R and P is a graded submodule of M with $(P:_RM) \cap S =\emptyset$. Then the following are equivalent:
\\
(i) P is a graded classical S-primary submodule of M.
\\
(ii) $S^{-1}P$ is a graded classical primary submodule of $S^{-1}M$ and there exists  $s \in S$ satisfying $(P:_Ms') \subseteq (P:_Ms)$ for all $s'\in S$.
\\
\\
\textbf{Proof.} (i) $\Rightarrow$ (ii): By Proposition \ref{prop3} and Lemma \ref{lem3}.
\\
(ii) $\Rightarrow$ (i): Let $x,y \in h(R)$ and $m \in h(M)$ with $xym \in P$. Then $\frac{x}{1}\frac{y}{1}\frac{m}{1} \in S^{-1}P$. Since $S^{-1}P$ is a graded classical primary submodule of M and M is a graded finitely generated, we get $\frac{x}{1}\frac{m}{1} \in S^{-1}P$ or $(\frac{y}{1})^n \frac{m}{1} \in S^{-1}P$ for some positive integer n. Then $s_1xm \in P$ or $s_2y^nm \in  P$ for some $s_1,s_2 \in S$. By assumption, there exists  $s \in S$ such that $(P:_Ms') \subseteq (P:_Ms)$ for all $s'\in S$. If $s_1xm \in P$, then $xm \in  (P:_M s_1) \subseteq (P:_Ms)$ and then $sxy \in P$. If $s_1y^nm \in P$, so $y^nm \in (P:_Ms_2) \subseteq (P:_Ms)$ and so $sy^nm \in P$ for some positive integer n. Therefore, P is a graded classical S-primary submodule of M. 
\end{theorem}

\begin{theorem}\label{thm8}
Let $S\subseteq h(R)$ be a multiplicative closed set of R and P is a graded submodule of R-module M. Then P is a graded classical S-primary submodule of M if and only if $(P:_Ms)$ is a graded classical primary submodule of M for some $s\in S$.
\\
\\
\textbf{Proof.} ($\Rightarrow$) Suppose that P is a graded classical S-primary submodule of M. Then there exists  $s\in S$ such that whenever $xym \in P$, where $x,y \in h(R)$ and $m \in h(M)$ and then $sxm \in P$ or $sy^n m \in P$ for some positive integer n. Now, let $abm \in (P:_Ms)$, where $a,b \in h(R)$. Then $sxym \in P$. As P is a graded classical S-primary submodule of M, we get $s^2xm \in P$ or $s^2y^nm \in P$. If $s^2xm \in P $, then we have nothing to show for it. Assume that $s^2xm \not\in P$ then $s^2y^n m \in P $ and then $sy^nm \in P$. Thus $y^nm \in (P:_Ms)$. Therefore, $(P:_M s)$ is a graded classical primary submodule of M for some $s\in S$.
\\
($\Leftarrow$) Suppose that $(P:_Ms)$ is a graded classical primary submodule of M for some $s\in S$. Let $xym \in P$ where $x,y \in h(R)$ and $m \in h(M)$. Then $xym \in (P:_Ms)$. Since $(P:_Ms)$ is a graded classical primary submodule of M, we get $xm \in (P:_Ms)$ or $y^nm \in (P:_Ms)$. which implies that either $sxm \in P$ or $sy^nm \in P$. Therefore, P is a graded classical S-primary submodule of M.
\end{theorem}

Note that in M a G-graded R-module over a G-graded ring R and I, J are graded submodules of M. Then $I \cap J$ is a graded submodule of M (see \cite{farzalipour2012union}).

\begin{proposition}\label{prop6}
Let $S \subseteq h(R)$ be a multiplicative closed set of Rand P is a graded classical S-primary submodule of R-module M such that $(P:_RM) \cap S =\emptyset$. If I is a graded submodule of M such that $(I:_RM) \cap S \neq \emptyset$, then $P\cap I$ is a graded classical S-primary submodule of M.
\\
\\
\textbf{Proof.} Since P and I is a graded submodules of M. then $P\cap I$ is a graded submodule of M. Clearly, $(P\cap I:_RM) \cap S = \emptyset$. Let $t \in (I:_RM) \cap S$ and $xym \in P \cap I$ where $x,y \in h(R)$ and $m \in h(m)$. Then $xym \in P$ and then $sxm \in P$ or $sy^nm\in P$ for some positive integer n and $s\in S$. Take $w = st \in S$. If $sxm \in P$ so $wxm \in P\cap I$ or if $sy^nm \in P$ so $wy^nm \in P \cap I$. Therefore, $P \cap I$ is a graded classical S-primary submodule of M. 
\end{proposition}

\begin{theorem}\label{thm9}
Let $S_1 , S_2 $ be a multiplicative closed subsets of rings $R_1 , R_2$ respectively and $P_1, P_2$ be a graded submodules of an $R_1$-module $M_1$ and $R_2$-module $M_2$ respectively. Suppose that $M = M_1 \times M_2$ as $R=R_1 \times R_2$-module, $S = S_1\times S_2$ and $P = P_1 \times P_2$. Then the following hold:
\\
(i) $P_1$ is a graded classical $S_1$-primary submodule of $M_1$ if and only if $P_1 \times M_2$ is a graded classical S-primary submodule of M.
\\
(ii) $P_2$ is a graded classical $S_2$-primary submodule of $M_2$ if and only if $M_1 \times P_2$ is a graded classical S-primary submodule of M.
\\
\\
\textbf{Proof.} (i) Let $(x,y)(a,b)(m_1,m_2)=(xam_1,ybm_2) \in P_1 \times M_2$ where $(x,y),(a,b) \in h(R)$ and $(m_1,m_2) \in h(M)$. Since $P_1$ is a graded classical S-primary submodule of $M_1$, then there exists  $s_1 \in S_1$  and then either $s_1xm_1 \in P$ or $s_1a^nm_1 \in P$ for some positive integer n. Let $s_2 \in S_2$, so $s_2ym_2 \in M_2$ and $s_2b^n M_2$. It follows that either $s(x,y)(m_1,m_2) \in P_1 \times M_2$ or $s(a^n,b^n)(m_1,m_2) \in P_1 \times M_2$ for some positive integer n. Therefore, $P_1\times M_2$ is a graded classical S-primary submodule of M. Conversely, take  $(x,y),(a,b) \in h(R)$ and $(m_1,m_2) \in h(M)$ and $(x,y)(a,b)(m_1,m_2)=(xam_1,ybm_2) \in P_1 \times M_2$ then $xam_1 \in P_1$. So there exists  $s=(s_1,s_2)\in S$ and so $s(x,y)(m_1,m_2)=(s_1xm_1,s_2ym_2) \in P_1 \times M_2$ or $s(a^n,b^n)(m_1,m_2)=(s_1a^nm_1,s_2b^nm_2) \in P_1 \times M_2$ for some positive integer n. Then $s_1xm_1 \in P_1$ or $s_1a^nm_1\in P_1$. Hence $P_1$ is a graded classical $S_1$-primary submodule of $M_1$.
\\
(ii) It is similar to (i).
\end{theorem}

\begin{corollary}\label{coro1}
Let $S_1 , S_2 $ be a multiplicative closed subsets of rings $R_1 , R_2$ respectively and $P_1, P_2$ be a graded submodules of an $R_1$-module $M_1$ and $R_2$-module $M_2$ respectively. Suppose that $M = M_1 \times M_2$ as $R=R_1 \times R_2$-module, $S = S_1\times S_2$ and $P = P_1 \times P_2$. Then the following are equivalent:
\\
(i) P is a graded classical S-primary submodule of M.
\\
(ii) $P_1 = M_1$ and $P_2$ is a graded classical $S_2$-primary submodule of $M_2$ or $P_2 = M_2$ and $P_1$ is a graded classical $S_1$-primary submodule of $M_1$.
\\
\\
\textbf{Proof.} It is follows by Theorem \ref{thm9}.
\end{corollary}

Note that if a G-graded ring R meets the ascending chain condition on graded ideals of R, it is termed a graded noetherian.

\begin{theorem}\label{thm10}
 Let $S\subseteq h(R)$ be multiplucative closed set of R and  P be  a graded submodule of M. If R is a noetherian ring, then every irreducible submodule of M is graded classical S-primary submodule of M.
 \\
 \\
 \textbf{Proof.} Let P be a irreducible submodule of M, $xym \in P$ anf $sxm \not \in P$, where $x,y \in h(R)$, $m \in h(M)$ and $s\in S$. Now, Let $(P:_R sym ) \subseteq (P:_R sy^2m) \subseteq ... \subseteq (P:_R sy^im) \subseteq ...$ be an ascending chain of ideals of R. Since R is a noetherian ring, there exists  a positive integer n such that $(P:_R sy^n m) = (P:_R sy^{n+i}m) $ for all $i \in \mathbf{N}$. We show that $P = (P+Rsy^nm) \cap (P+Rsxm)$. Let $t \in (P+Rsy^nm) \cap (P+Rsxm)$.  Then $t = p_1 + r_1sym = p_2 + r_2sxm$, where $p_1 ,p_2 \in P$ and $r_1 , r_2 \in h(R)$ and hence $yt = yp_1 + r_1sy^{n+1}m = yp_2 + yr_2sxm \in  P$. Since $yp_1, yp_2, yr_2sxm \in P$, we have $r_1 sy^{n+1}m \in P$. So $p_1 \in (P:_R sy^n m) = (P:_R sy^{n+1}m) $, i.e, $p_1sy^nm \in P$ and thus $t = p_1 + r_1sy^nm \in P$. Hence we have $(P + Rsy^nm) \cap (P + Rsxm) \subseteq P$. The inclusion is opposite side is self-evident, therefore $P= (P + Rsy^nm) \cap (P + Rsxm)$. Since P is an irreducible submodule and $sxm \not\in P$, we have $P = (P + Rsy^nm)$. Hence $sy^nm \in P$. Therefore P is a  graded classical S-primary submodule of M.  
\end{theorem}

\bibliographystyle{amsplain}

\begin{thebibliography}{99}

\bibitem{abu2019graded} Rashid Abu-dawwas and Malik Bataineh. \textit{“Graded r-ideals”}.  In: Iranian Journal of
Mathematical Sciences and Informatics 14.2 (2019), pp. 1–8.

\bibitem{atani2006graded} S Ebrahimi Atani.  \textit{”On graded prime submodules”}. In: Chiang Mai J. Sci 33.1 (2006),
pp. 3–7.

\bibitem{atani2007graded} Shahabaddin E Atani and Unsal Tekir. \textit{”On the graded primary avoidance theorem”}. In: Chiang Mai J. Sci 34.2 (2007), pp. 161–164.

\bibitem{escoriza1998multiplication} J Escoriza and B Torrecillas. \textit{“Multiplication objects in commutative Grothendieck
categories”}.   In: Communications in algebra 26.6 (1998), pp. 1867–1883.

\bibitem{farzalipour2012union} Farkhonde Farzalipour and Peyman Ghiasvand.  \textit{ “On the union of graded prime submodules”}. In: Thai journal of mathematics 9.1 (2012), pp. 49–55.


\bibitem{gilmer1992multiplicative} Robert Gilmer.  \textit{ “Multiplicative ideal theory”}.  In: Queen’s papers in Pure and Appl. Math. 90 (1992).


\bibitem{gordon1982graded} Robert Gordon and Edward L Green.  \textit{ “Graded artin algebras”}. In: Journal of Algebra 76.1 (1982), pp. 111–137


\bibitem{nastasescu2004methods} Constantin Nastasescu, Freddy Van Oystaeyen, and Freddy MJ Van Oystaeyen.  \textit{”Methods of graded rings”}. 1836. Springer Science and Business Media, 2004.


\bibitem{refai2000graded} Mashhoor Refai. \textit{“Graded radicals and graded prime spectra”}. In: Far East journal of mathematical sciences (2000), pp. 59–73.  


\bibitem{refai2004graded} Mashhoor Refai and Khaldoun Al-Zoubi. \textit{ “On graded primary ideals”}. In: Turkish journal of mathematics 28.3 (2004), pp. 217–230. 


\bibitem{al2022generalizations} Tamem Al-Shorman, Malik Bataineh, and Rashid Abu-Dawwas. \textit{“Generalizations of Graded S-Primary Ideals”}. In: arXiv preprint arXiv:2203.04092 (2022).  


\bibitem{al2016graded} Khaldoun Al-Zoubi. \textit{“The graded primary radical of a graded submodules”}. In: An. Stiint. Univ. Al. I. Cuza Iasi. Mat.(NS) 1 (2016), pp. 395–402.


\bibitem{al2021graded} Khaldoun Al-Zoubi and Farah Ababneh. \textit{“ON GRADED S-PRIMARY SUBMODULES OF GRADED MODULES OVER GRADED COMMUTATIVE RINGS”}. In:
(2021). 


\bibitem{kmm2021graded} Khaldoun Al-Zoubi, Mohammed Ali, and Mohammad Alkhatib. \textit{“On graded classical S-2-absorbing submodules”}. In: (2021).


\bibitem{al2017graded} Khaldoun Al-Zoubi and Mohammed Al-Dolat. \textit{“On graded classical primary submodules”}. In: Advances in pure and applied mathematics 7.2 (2016), pp. 93–96.
11.






\end{thebibliography}

\end{document}